\documentclass{article}
\usepackage{setspace}
\usepackage{amsmath}
\usepackage{amssymb}
\usepackage{hyphenat}
\usepackage{amsfonts}
\usepackage{url}
\usepackage[T1]{fontenc}
\usepackage{ae}
\usepackage{aecompl}
\newtheorem{theorem}{Theorem}
\newtheorem{definition}{Definition}
\newtheorem{example}{Example}
\newtheorem{remark}{Remark}
\title{An epsilon-delta characterization\\of a certain TTE computability notion}
\author{Dimiter Skordev}
\date{}
\begin{document}
\maketitle

\begin{abstract}
The TTE computability notion in effective metric spaces is usually defined by using Cauchy representations. Under some weak assumptions, we characterize this notion in a way which avoids using the representations.
\end{abstract}

\section{Introduction}

The widely used TTE approach to computability of functions in $\mathbb{R}$ and in other non{\hyp}denumerable domains uses computable transformations of infinitistic names of the argument values into the same kind of names of the function's values. In addition, general quantifiers over the names of the argument values are used in this approach. In some cases, the use of such infinitistic names can be avoided (cf., for instance, \cite[Theorem 3]{skor}, where certain TTE computable real functions are characterized in the spirit of the notion from \cite{TZ} of being uniformly in a class of total functions in $\mathbb{N}$). However, as far as we know, no characterization of this kind is presented yet for the general TTE computability notion for real functions which is based on Cauchy representations of the real numbers. 

In Section \ref{comp}, which is the main part of this paper, effective metric spaces and their Cauchy representations in the sense of \cite{weih} will be considered. Under some weak assumptions, the corresponding TTE computability notion will be characterized in a way which avoids using the representations and can be regarded in some sense as an epsilon{\hyp}delta approach. To give the characterization in question, we will initially (in Section \ref{cont}) consider arbitrary pairs of metric spaces and, for any partial mapping $\theta$ of the first space into the second one, we define the notion of an $A,B${\hyp}approximation system for $\theta$, whenever $A$ and $B$ are dense subsets of the first and the second space, respectively (the partial mappings $\theta$ having approximation systems will turn out to be exactly the continuous ones).

\section{Approximation systems}\label{cont}

Throughout the paper, metric spaces $(X,d)$ and $(Y,e)$ are supposed to be given, as well as a function $\theta:D\to Y$, where $D\subseteq X$.

\begin{definition}\label{apprsyst}
Let $A$ and $B$ be dense subsets of the metric spaces $(X,d)$ and $(Y,e)$, respectively. A subset $S$ of $A\times\mathbb{N}\times B\times\mathbb{N}$ will be called an {\em $A,B${\hyp}approximation system} for $\theta$ if the following two conditions are satisfied for any $\xi\in D$:
\begin{enumerate}
\item Whenever $(a,m,b,n)\in S$ and $d(a,\xi)<\frac{1}{m+1}$ the inequality \mbox{$e(b,\theta(\xi))<\frac{1}{n+1}$} holds.
\item For any $n\in\mathbb{N}$, there exists $m\in\mathbb{N}$ such that whenever $a\in A$ and \mbox{$d(a,\xi)<\frac{1}{m+1}$} an element $b$ exists such that $(a,m,b,n)\in S$. 
\end{enumerate}
\end{definition}

\begin{example}\label{division}
{\em Let
\begin{gather*}
X=\mathbb{R}^2,\ \ d((\xi_1,\xi_2),(\xi_1^\prime,\xi_2^\prime))=\max(|\xi_1{\hyp}\xi_1^\prime|,|\xi_2{\hyp}\xi_2^\prime|),\\Y=\mathbb{R},\ \ e(\eta,\eta^\prime)=|\eta{\hyp}\eta^\prime|,\ \ D=\mathbb{R}\times(\mathbb{R}\setminus\{0\}),\ \ \theta(\xi_1,\xi_2)=\frac{\xi_1}{\xi_2}.
\end{gather*}
Let $S$ be the set of all elements $((a_1,a_2),m,b,n)$ of $\mathbb{Q}^2\times\mathbb{N}\times\mathbb{Q}\times\mathbb{N}$ such that $a_2b=a_1$ and $(m+1)|a_2|\ge 1+(n+1)(|b|+1)$.
Then $S$ is a $\mathbb{Q}^2,\mathbb{Q}${\hyp}approximation system for $\theta$. To prove this, suppose $(\xi_1,\xi_2)\in D$. If $((a_1,a_2),m,b,n)\in S$ and $d((a_1,a_2),(\xi_1,\xi_2))<\frac{1}{m+1}$ then $$a_2\ne 0,\ \ b=\frac{a_1}{a_2},\ \ |\xi_2|>|a_2|-\frac{1}{m+1}\ge\frac{(n+1)(|b|+1)}{m+1},$$
\vspace{-5mm}
\begin{multline*}
|b-\theta(\xi_1,\xi_2)|=\left|\frac{a_1}{a_2}-\frac{\xi_1}{\xi_2}\right|=\frac{|a_1\xi_2-a_2\xi_1|}{|a_2||\xi_2|}\\
\le\frac{|a_1||\xi_2-a_2|+|a_2||a_1-\xi_1|}{|a_2|(n+1)(|b|+1)}(m+1)<
\frac{|a_1|+|a_2|}{|a_2|(n+1)(|b|+1)}=\frac{1}{n+1}.
\end{multline*}
Suppose now an arbitrary natural number $n$ is given. Let $m$ be a natural number such that $(m+1)\xi_2^2\ge 3|\xi_2|+1+(n+1)(|\xi_1|+|\xi_2|+2)$, and let $a_1,a_2$ be rational numbers satisfying $d((a_1,a_2),(\xi_1,\xi_2))<\frac{1}{m+1}$. Obviously $(m+1)\xi_2^2>|\xi_2|$, hence $|\xi_2|>\frac{1}{m+1}$ and therefore $a_2\ne 0$. We will show that $((a_1,a_2),m,b,n)\in S$ for $b=a_1/a_2$. Indeed, then $a_2b=a_1$ and
\begin{multline*}
1\!+\!(n\!+\!1)(|b|\!+\!1)=\frac{|a_2|\!+\!(n\!+\!1)(|a_1|\!+\!|a_2|)}{|a_2|}<\frac{|\xi_2|\!+\!1\!+\!(n\!+\!1)(|\xi_1|\!+\!|\xi_2|\!+\!2)}{|a_2|}\\
\le\frac{(m+1)\xi_2^2-2|\xi_2|}{|a_2|}<\frac{m+1}{|a_2|}\left(|\xi_2|-\frac{1}{m+1}\right)^2<\frac{m+1}{|a_2|}|a_2|^2=(m+1)|a_2|.
\end{multline*}
}
\end{example}

\begin{example}\label{cosine}{\em Let $(X,d)=(Y,e)$, where $(Y,e)$ is the same as in Example \ref{division}, $D$ be the set $\mathbb{R}$, and $\theta$ be defined by $\theta(\xi)=\cos\xi$. Making use of the inequality \mbox{$|\cos a-\cos\xi|\le|a-\xi|$} and the equality
$$\cos a=\sum_{i=0}^\infty(-1)^i\frac{a^{2i}}{(2i)!},$$
one can conclude that, for any $\xi\in\mathbb{R}$, any $a\in\mathbb{Q}$ and any $k\in\mathbb{N}$ such that $(2k+1)(2k+2)\ge a^2$, the inequality
$$\left|\sigma_k(a)-\cos\xi\right|\le\frac{a^{2k}}{2(2k)!}+|a-\xi|$$
holds, where
$$\sigma_k(a)=(-1)^k\frac{a^{2k}}{2(2k)!}+\sum_{i<k}(-1)^i\frac{a^{2i}}{(2i)!}.$$
For any $k\in\mathbb{N}$, let $S_k$ be the set of all $(a,m,b,n)\in\mathbb{Q}\times\mathbb{N}\times\mathbb{Q}\times\mathbb{N}$ such that
$$a^2\le(2k+1)(2k+2),\ \ |b-\sigma_k(a)|+\frac{a^{2k}}{2(2k)!}+\frac{1}{m+1}\le\frac{1}{n+1}.$$
Then the union of the sets $S_0,S_1,S_2,\ldots$ is a $\mathbb{Q},\mathbb{Q}${\hyp}approximation system for $\theta$. To prove this, suppose $\xi$ is an arbitrary real number. If $(a,m,b,n)\in S_k$ for some $k\in\mathbb{N}$ and $d(a,\xi)<\frac{1}{m+1}$ then
\begin{multline*}
e(b,\theta(\xi))\le|b-\sigma_k(a)|+|\sigma_k(a)-\cos\xi|
\le|b-\sigma_k(a)|+\frac{a^{2k}}{2(2k)!}+|a-\xi|\\
<|b-\sigma_k(a)|+\frac{a^{2k}}{2(2k)!}+\frac{1}{m+1}\le\frac{1}{n+1}.
\end{multline*}
Suppose now a natural number $n$ is given. Let $k$ be a natural number such that
$$(2k+1)(2k+2)\ge(|\xi|+1)^2,\ \ \frac{(|\xi|+1)^{2k}}{2(2k)!}<\frac{1}{n+1},$$
and let $m\in\mathbb{N}$ satisfy the inequality
$$\frac{(|\xi|+1)^{2k}}{2(2k)!}+\frac{1}{m+1}\le\frac{1}{n+1}.$$
If $a$ is any rational number such that $d(a,\xi)<\frac{1}{m+1}$ then $|a|<|\xi|+1$ and therefore
$$a^2<(|\xi|+1)^2\le(2k+1)(2k+2),\ \ \frac{a^{2k}}{2(2k)!}+\frac{1}{m+1}<\frac{(|\xi|+1)^{2k}}{2(2k)!}+\frac{1}{m+1}\le\frac{1}{n+1},$$
hence $(a,m,n,\sigma_k(a))\in S_k$.
}\end{example}

\begin{theorem}\label{continuous}
Let $A$ and $B$ be dense subsets of the metric spaces $(X,d)$ and $(Y,e)$, respectively. An $A,B${\hyp}approximation system for $\theta$ exists iff $\theta$ is continuous, and if $\theta$ is continuous then the set of all $(a,m,b,n)\in A\times\mathbb{N}\times B\times\mathbb{N}$ satisfying the condition
\begin{equation}\label{max}
\forall\xi\in D\left(d(a,\xi)<\frac{1}{m+1}\Rightarrow e(b,\theta(\xi))<\frac{1}{n+1}\right)
\end{equation}
is an $A,B${\hyp}approximation system for $\theta$.
\end{theorem}

{\em Proof.} For proving the implication from left to right in the first part of the conclusion, suppose $S$ is an $A,B${\hyp}approximation system for $\theta$. Let $\xi\in D$, and $\varepsilon$ be a positive number. We choose a natural number $n$ satisfying the inequality $\frac{2}{n+1}<\varepsilon$. By condition~2 of Definition~\ref{apprsyst}, a natural number $m$ exists such that whenever $a\in A$ and \mbox{$d(a,\xi)<\frac{1}{m+1}$} the quadruple $(a,m,b,n)$ belongs to $S$ for some $b$. Let $m$ be such a natural number, and let $\xi^\prime$ be any element of $D$ satisfying the inequality $d(\xi^\prime,\xi)<\frac{1}{2m+2}$. After choosing an element $a$ of $A$, satisfying the inequality \mbox{$d(a,\xi^\prime)<\frac{1}{2m+2}$,} we will have both inequalities $d(a,\xi^\prime)<\frac{1}{m+1}$ and $d(a,\xi)<\frac{1}{m+1}$. By the second of them, $(a,m,b,n)\in S$ for some $b$. The two inequalities and condition 1 of Definition~\ref{apprsyst} imply the inequalities $e(b,\theta(\xi^\prime))<\frac{1}{n+1}$ and $e(b,\theta(\xi))<\frac{1}{n+1}$, hence \mbox{$e(\theta(\xi^\prime),\theta(\xi))<\frac{2}{n+1}<\varepsilon$.} Thus the continuity of $\theta$ is established. For proving the rest of the conclusion, suppose now $\theta$ is continuous. Let $S$ be the set of all $(a,m,b,n)\in A\times\mathbb{N}\times B\times\mathbb{N}$ satisfying the condition (\ref{max}). We will prove that $S$ is an $A,B${\hyp}approximation system for~$\theta$. Let $\xi$ be an arbitrary element of~$D$. Condition 1 of Definition~\ref{apprsyst} follows immediately from the definition of~$S$. To check condition 2, suppose $n$ is a natural number. By the continuity of $\theta$, a positive number $\delta$ exists such that $e(\theta(\xi^\prime),\theta(\xi))<\frac{1}{2n+2}$ for all $\xi^\prime$ in~$D$ satisfying the inequality $d(\xi^\prime,\xi)<\delta$. We choose a natural number $m$ with $\frac{2}{m+1}<\delta$ and an element $b$ of $B$ such that $e(b,\theta(\xi))<\frac{1}{2n+2}$. Consider now any $a$ in~$A$ satisfying the inequality $d(a,\xi)<\frac{1}{m+1}$. We will show that $(a,m,b,n)\in S$. To do this, suppose $\xi^\prime$ is any element of $D$ satisfying the inequality $d(a,\xi^\prime)<\frac{1}{m+1}$. Then 
$$d(\xi^\prime,\xi)\le d(\xi^\prime,a)+d(a,\xi)<\frac{2}{m+1}<\delta$$
and consequently $e(\theta(\xi^\prime),\theta(\xi))<\frac{1}{2n+2}$. Since 

$$e(b,\theta(\xi^\prime))\le e(b,\theta(\xi))+e(\theta(\xi),\theta(\xi^\prime)),$$
we see that $e(b,\theta(\xi^\prime))<\frac{1}{n+1}$.$\ \ _\Box$

\begin{remark}\label{maximal}
{\em If $\theta$ is continuous then the set of all $(a,m,b,n)\in A\times\mathbb{N}\times B\times\mathbb{N}$ satisfying the condition (\ref{max}) is obviously the maximal $A,B${\hyp}approximation system for $\theta$ -- it contains as subsets all other ones. Let us note that, as seen from the above proof, condition 2 of Definition \ref{apprsyst} is satisfied for this set in a stronger form, namely $b$ does not depend on the choice of $a$.}
\end{remark}

\begin{example}\label{division-max}
{\em Let $X,\,d,\,Y,\,e,\,\theta$ be as in Example \ref{division}. Then the maximal $\mathbb{Q}^2,\mathbb{Q}${\hyp}approximation system for $\theta$ consists of all elements $((a_1,a_2),m,b,n)$ of the set \mbox{$\mathbb{Q}^2\times\mathbb{N}\times\mathbb{Q}\times\mathbb{N}$} such that $(m+1)|a_2|>1$ and the four numbers of the form
$$\frac{(m+1)a_1\pm 1}{(m+1)a_2\pm 1}$$
belong to the closed interval $\left[\,b-\frac{1}{n+1},\,b+\frac{1}{n+1}\,\right]$.}
\end{example}

\section{Characterization of the TTE computability \\by means of approximation systems}\label{comp}

In order to be able to discuss computability questions about the function~$\theta$, we will now suppose some additional structure on the metric spaces $(X,d)$ and $(Y,e)$ turning them into effective metric spaces $\mathbf{X}=(X,d,A,\alpha)$, \mbox{$\mathbf{Y}=(Y,e,B,\beta)$} in the sense of  \cite[Definition~8.1.2]{weih} (the same definition introduces also the notion of Cauchy representation associated with an effective metric space). Namely, $A$ and $B$ must be countable dense subsets of $(X,d)$ and $(Y,e)$, respectively, and $\alpha$, $\beta$ must be notations of them in the sense of \mbox{\cite[Definition~2.3.1]{weih}} (unlike the definition given in \cite{hemm}, Definition~8.1.2 in \cite{weih} does not require completeness of the considered metric space, imposes no effectivity requirement on its metric and allows arbitrary sets of finite words over a finite alphabet to be domains of the namings). Definition 3.1.3 of \cite{weih} allows using the notion of $(\delta_\mathbf{X},\delta_\mathbf{Y})${\hyp}computability of the function $\theta$, where $\delta_\mathbf{X}$ and $\delta_\mathbf{Y}$ are the Cauchy representations associated with $\mathbf{X}$ and $\mathbf{Y}$. For the sake of technical convenience, the notations $\alpha$ and $\beta$ will be regarded as mappings of certain subsets of $\mathbb{N}$ onto $A$ and onto $B$, respectively. Then any Cauchy name of an element $\xi$ of $X$ can be identified with a total one{\hyp}argument function $f$ in $\mathbb{N}$ such that $f(i)\in\mathrm{dom}(\alpha)$ for any $i\in\mathbb{N}$, $d(\alpha(f(i)),\alpha(f(k)))\le 2^{-i}$ for $i<k$, and $\xi=\lim_{i\to\infty}\alpha(f(i))$ (similarly for the Cauchy names of the elements of~$Y$). The $(\delta_\mathbf{X},\delta_\mathbf{Y})${\hyp}computability of the function~$\theta$ will be equivalent to the existence of a recursive operator $T$ such that $T$ transforms unary partial functions in $\mathbb{N}$ into unary partial functions in $\mathbb{N}$ and, whenever $f$ is a Cauchy name of an element $\xi$ of $D$, $T(f)$ is a Cauchy name of the element $\theta(\xi)$ of~$Y$.

In the sequel, it would be convenient to use some other names of the elements of $X$ and $Y$. In the situation described above, an {\em ordinary name} of an element $\xi$ of $X$ will be, by definition, a total one{\hyp}argument function $f$ in $\mathbb{N}$ such that \mbox{$f(i)\in\mathrm{dom}(\alpha)$} and $d(\alpha(f(i)),\xi)<\frac{1}{i+1}$ for any $i\in\mathbb{N}$ (similarly for the elements of $Y$). Since there exist a primitive recursive operator transforming the ordinary names of any element into Cauchy names of the same element, as well as a primitive recursive operator performing a transformation in the opposite direction, the $(\delta_\mathbf{X},\delta_\mathbf{Y})${\hyp}computability of the function $\theta$ will be also equivalent to the existence of a recursive operator $T$ such that $T$ transforms unary partial functions in $\mathbb{N}$ into unary partial functions in $\mathbb{N}$ and, whenever $f$ is an ordinary name of an element $\xi$ of $D$, $T(f)$ is an ordinary name of the element $\theta(\xi)$ of~$Y$.

\begin{definition}\label{enumerable}
Let $S\subseteq A\times\mathbb{N}\times B\times\mathbb{N}$. The set $S$ will be called {\em recursively enumerable} with respect to $\alpha$ and $\beta$ if the set
$$\hat{S}=\{(i,m,j,n)\,|\,i\in\mathrm{dom}(\alpha),\,j\in\mathrm{dom}(\beta),\,(\alpha(i),m,\beta(j),n)\in S\}$$
is a recursively enumerable subset of $\mathbb{N}^4$.
\end{definition}

\begin{example}\label{example_for_enumerable}
{\em Let $A=\mathbb{Q}^N$, $B=\mathbb{Q}$, and let the notations $\alpha$ of $A$ and $\beta$ of $B$ be representable in the following form, where $\rho_1,\sigma_1,\tau_1,\ldots,\rho_N,\sigma_N,\tau_N,\rho,\sigma,\tau$ are some partial recursive functions:
$$\alpha(i)=\left(\frac{\rho_1(i)-\sigma_1(i)}{\tau_1(i)+1},\ldots,\frac{\rho_N(i)-\sigma_N(i)}{\tau_N(i)+1}\right),\ \ \beta(j)=\frac{\rho(j)-\sigma(j)}{\tau(j)+1}$$
($\mathrm{dom}(\alpha)$ and $\mathrm{dom}(\beta)$ are supposed to be the intersection of the domains of $\rho_1,\sigma_1,\tau_1,\ldots,\rho_N,\sigma_N,\tau_N$ and the intersection of the domains of $\rho$, $\sigma$, $\tau$, respectively). Let $S\subseteq A\times\mathbb{N}\times B\times\mathbb{N}$, and let $S^\dag$ be the set of all elements 
$(r_1,s_1,t_1,\ldots,r_N,s_N,t_N,m,r,s,t,n)$ of $\mathbb{N}^{3N+5}$ such that
$$\left(\left(\frac{r_1-s_1}{t_1+1},\ldots,\frac{r_N-s_N}{t_N+1}\right),m,\frac{r-s}{t+1},n\right)\in S.$$
The set $S$ is recursively enumerable with respect to $\alpha$ and $\beta$ iff $S^\dag$ is a recursively enumerable subset of $\mathbb{N}^{3N+5}$.}
\end{example}

\begin{theorem}\label{ex_reas->comp}
Let $S$ be an $A,B${\hyp}approximation system for the function $\theta$, and let $S$ be recursively enumerable with respect to $\alpha$ and $\beta$. Then $\theta$ is $(\delta_\mathbf{X},\delta_\mathbf{Y})${\hyp}computable.
\end{theorem}

{\em Proof.} By the recursive enumerability of $S$ with respect to $\alpha$ and $\beta$, a 5{\hyp}argument primitive recursion function $\chi$ can be found such that, for all $i,m,j,n\in\mathbb{N}$, the equivalence $(i,m,j,n)\in\hat{S}\Leftrightarrow\exists s(\chi(i,m,j,n,s)=0)$ holds. Let $\pi_1,\pi_2,\pi_3$ be unary primitive recursive functions such that
$$\{(\pi_1(k),\pi_2(k),\pi_3(k))\,|\,k\in\mathbb{N}\}=\mathbb{N}^3.$$ Let us define recursive operators $K$ and $T$ as follows:
\begin{align*}
K(f)(n)=&\ \mu k[\,\chi(f(\pi_1(k)),\pi_1(k),\pi_2(k),n,\pi_3(k))=0\,],\\
T(f)(n)=&\ \pi_2(K(f)(n)).
\end{align*}
We will prove the $(\delta_\mathbf{X},\delta_\mathbf{Y})${\hyp}computability of $\theta$ by showing that, whenever \mbox{$\xi\in D$} and $f$ is an ordinary name of $\xi$, the function $T(f)$ is an ordinary name of $\theta(\xi)$. Let $f$ be an ordinary name of an element $\xi$ of $D$, and let $n\in\mathbb{N}$. Making use of condition 2 of Definition~\ref{apprsyst}, we choose a natural number $m$ such that, whenever $a\in A$ and $d(a,\xi)<\frac{1}{m+1}$, the quadruple $(a,m,b,n)$ belongs to $S$ for some $b\in B$. Since \mbox{$d(\alpha(f(m)),\xi)<\frac{1}{m+1}$,} there exists \mbox{$j\in\mathrm{dom}(\beta)$} such that $(\alpha(f(m)),m,\beta(j),n)\in S$, hence $(f(m),m,j,n)\in\hat{S}$ and therefore \mbox{$\chi(f(m),m,j,n,s)=0$} for some $s\in\mathbb{N}$. It follows from here that $n$ belongs to $\mathrm{dom}(K(f))$, and consequently $n$ belongs also to $\mathrm{dom}(T(f))$. After setting $K(f)(n)=k$, $\pi_1(k)=l$, we will have the equality
$$\chi(f(l),l,T(f)(n),n,\pi_3(k))=0,$$
and it shows that $(f(l),l,T(f)(n),n)\in\hat{S}$, i.e.\@ $f(l)$ and $T(f)(n)$ belong to $\mathrm{dom}(\alpha)$ and $\mathrm{dom}(\beta)$, respectively, and $(\alpha(f(l)),l,\beta(T(f)(n)),n)$ belongs to $S$. From here, making use of condition 1 of Definition~\ref{apprsyst}, we conclude that
$$e(\beta(T(f)(n)),\theta(\xi))<\frac{1}{n+1}.$$
Since this reasoning was done for an arbitrary natural number $n$, we thus proved that $T(f)$ is really an ordinary name of~$\theta(\xi).\ \ _\square$

\vskip3mm
Under some weak additional assumptions, a converse of the above theorem also holds (the assumptions in question are surely satisfied if the effective metric spaces $\mathbf{X}$ and $\mathbf{Y}$ are semi{\hyp}computable in the sense introduced on page 239 \mbox{of \cite{weih}).}

\vspace{2mm}
\begin{theorem}\label{comp->ex_re_as}
\hspace{-1.7mm}{\bf.} Let the sets
\begin{align}
\left\{(p,q,r)\left|\ p,q\in\mathrm{dom}(\alpha),\ r\in\mathbb{N},\ d(\alpha(p),\alpha(q))<\frac{1}{2r+2}\right.\right\}&,\label{a}\\ \left\{(p,q,r)\left|\ p,q\in\mathrm{dom}(\beta),\ r\in\mathbb{N},\ e(\beta(p),\beta(q))<\frac{1}{2r+2}\right.\right\}&\label{b}
\end{align}
be recursively enumerable subsets of $\mathbb{N}^3$, and let the function $\theta$ be $(\delta_\mathbf{X},\delta_\mathbf{Y})${\hyp}computable. Then there exists an $A,B${\hyp}approximation system $S$ for $\theta$ such that $S$ is recursively enumerable with respect to $\alpha$ and $\beta$.
\end{theorem}

{\em Proof.} Let $T$ be a recursive operator such that $T$ transforms unary partial functions in $\mathbb{N}$ into unary partial functions in $\mathbb{N}$ and, whenever $f$ is an ordinary name of an element $\xi$ of $D$, the function $T(f)$ is an ordinary name of $\theta(\xi)$. Let $S$ be the set of all quadruples $(a,m,b,n)\in A\times\mathbb{N}\times B\times\mathbb{N}$ such that, for some $l\in\mathbb{N}$ satisfying the inequality $2l+1\le m$, a function $g$ from $\{0,1,\ldots,l\}$ into $\mathrm{dom}(\alpha)$ exists which satisfies the following conditions:
\begin{gather}
d(\alpha(g(k)),a)<\frac{1}{2k+2},\ \ k=0,1,\ldots,l,\label{d}\\
2n+1\in\mathrm{dom}(T(g)),\label{in}\\
T(g)(2n+1)\in\mathrm{dom}(\beta),\ \ e(\beta(T(g)(2n+1)),b)<\frac{1}{2n+2}\label{me}.
\end{gather}
The set $S$ is recursively enumerable with respect to $\alpha$ and $\beta$ due to the recursiveness of the operator $T$ and the recursive enumerability of the sets (\ref{a}) and (\ref{b}). We will show that $S$ is an $A,B${\hyp}approximation system for~$\theta$. 
Let $\xi\in D$. To verify condition 1 of Definition \ref{apprsyst}, suppose that $(a,m,b,n)\in S$ and $d(a,\xi)<\frac{1}{m+1}$. Let $l$ and $g$ be a natural number and a function with the properties from the above definition of the set $S$. For any $k\in\{0,1,\ldots,l\}$, we have
\begin{equation*}
d(\alpha(g(k)),\xi)\le d(\alpha(g(k)),a)+d(a,\xi)<\frac{1}{2k+2}+\frac{1}{2l+2}\le\frac{1}{k+1}.
\end{equation*}
It follows from here that $g$ can be extended to some ordinary name~$f$ of $\xi$, and then $T(f)$ will be an ordinary name of $\theta(\xi)$. By condition (\ref{in}) and the continuity of the operator~$T$, the equality
\begin{equation}\label{eq-restr}
T(f)(2n+1)=T(g)(2n+1)
\end{equation}
holds. Therefore, making use also of condition (\ref{me}), we have
\begin{multline*}
e(b,\theta(\xi))\le e(b,\beta(T(g)(2n+1)))+e(\beta(T(f)(2n+1),\theta(\xi))\\<\frac{1}{2n+2}+\frac{1}{2n+2}=\frac{1}{n+1}.
\end{multline*}
To verify condition 2, suppose a natural number~$n$ is given. Let \mbox{$f:\mathbb{N}\to\mathrm{dom}(\alpha)$} be such that $d(\alpha(f(k)),\xi)<\frac{1}{4k+4}$ for any $k\in\mathbb{N}$. Clearly $f$ is an ordinary name of~$\xi$, hence the function $T(f)$ is an ordinary name of $\theta(\xi)$, and therefore $T(f)$ is total and all its values belong to $\mathrm{dom}(\beta)$. Let $b=\beta(T(f)(2n+1))$. By the continuity of the operator~$T$, a natural number $l$ exists such that $g=f\upharpoonright\{0,1,\ldots,l\}$ satisfies condition~(\ref{in}) and the equality (\ref{eq-restr}), hence \mbox{$b=\beta(T(g)(2n+1))$.} Let $m=4l+3$, and let $a$ be an arbitrary element of~$A$ satisfying the inequality $d(a,\xi)<\frac{1}{m+1}$. Then  the inequalities~(\ref{d}) hold, because
$$d(\alpha(g(k)),a)\le d(\alpha(g(k)),\xi)+d(\xi,a)<\frac{1}{4k+4}+\frac{1}{4l+4}\le\frac{1}{2k+2}$$
for any $k\in\{0,1,\ldots,l\}$. Since obviously $m>2l+1$, and condition (\ref{me}) is trivially satisfied, we see that $(a,m,b,n)\in S.\ \ _\square$

\begin{remark}
{\em The set $S$ constructed in the above proof satisfies condition 2 of Definition~\ref{apprsyst} in the stronger form mentioned in Remark \ref{maximal} (the element $b$ does not depend on the choice of $a$). However, this set is not necessarily the maximal approximation system for $\theta$. Moreover, as seen from the next example, the assumptions of Theorem \ref{comp->ex_re_as} do not imply that the maximal approximation system for $\theta$ is necessarily recursively enumerable.} 
\end{remark}

\begin{example}{\em
Let $X=A=Y=B=\mathbb{N}$, the metrics $d$ and $e$ coincide with the usual metric in $\mathbb{N}$ (i.e. the absolute value of the difference), and $\alpha=\beta=\mathrm{id}_\mathbb{N}$. Let $D$ be such that $\mathbb{N}\setminus D$ is not recursively enumerable, $\theta$ be the restriction of the constant 0 to $D$, and $S$ be the maximal $A,B${\hyp}approximation system for $\theta$. Then, for any $a,m,b,n\in\mathbb{N}$, the condition (\ref{max}) is equivalent to the implication $a\in D\Rightarrow b=0$. Therefore $(a,0,1,0)\in S\Leftrightarrow a\not\in D$ for any $a\in\mathbb{N}$, hence $S$ is not recursively enumerable.
}\end{example}

\section{An application to TTE computability \\of semialgebraic real functions}
A real function will be called {\em strongly semialgebraic} if it is semicomputable in the sense of \cite[Section 4]{TZ}, i.e. if it is definable without parameters in the ordered field~$\mathbb{R}$. Equivalently (in virtue of quantifier elimination), a real function is strongly semialgebraic iff its graph can be defined by means of a Boolean combination of polynomial inequalities with integer coefficients.

\begin{theorem}\label{semialg}
Any continuous strongly semialgebraic real function is computable in the sense of \cite[Section 4.3]{weih}.
\end{theorem}

{\em Proof.} Let $\theta:D\to\mathbb{R}$, where $D\subseteq\mathbb{R}^N$, be continuous and strongly semialgebraic. We consider effective metric spaces $\mathbf{X}=(X,d,A,\alpha)$, $\mathbf{Y}=(Y,e,B,\beta)$, where $X=\mathbb{R}^N$, $Y=\mathbb{R}$, $e$ is the same as in Example \ref{division}, $d$ is defined similarly to the metric $d$ in that example, and $A,\alpha,B,\beta$ are as in Example \ref{example_for_enumerable}. Let $S$ be the set of all $(a,m,b,n)\in A\times\mathbb{N}\times B\times\mathbb{N}$ satisfying the condition (\ref{max}). By Theorem \ref{continuous}, $S$ is an $A,B${\hyp}approximation system for $\theta$. With $a=(a_1,\ldots,a_N)$, where $a_1,\ldots,a_N\in\mathbb{Q}$, the condition (\ref{max}) can be written in the form
\begin{multline*}
\forall\xi_1\in\mathbb{R}\,\ldots\,\forall\xi_N\in\mathbb{R}\,\forall\eta\in\mathbb{R}\,\bigg((\xi_1,\ldots,\xi_N,\eta)\in\varTheta\\
\&\ |a_1-\xi_1|<\frac{1}{m+1}\ \&\ \ldots\ \&\ |a_N-\xi_N|<\frac{1}{m+1}\,\Rightarrow\,|b-\eta|<\frac{1}{n+1}\bigg),
\end{multline*}
where $\varTheta$ is the graph of $\theta$. Making use of quantifier elimination, we can characterize $S$ by means of a Boolean combination of inequalities of the form 
$$P\left(a_1,\ldots,a_N,b,\frac{1}{m+1},\frac{1}{n+1}\right)>0,$$
where all $P$ are polynomials with integer coefficients. Therefore (by Example~\ref{example_for_enumerable}) $S$ is recursively enumerable with respect to $\alpha$ and $\beta$, hence Theorem \ref{ex_reas->comp} yields the $(\delta_\mathbf{X},\delta_\mathbf{Y})${\hyp}computability of $\theta$, which, of course, implies the computability of~$\theta$ in the sense of \cite[Section 4.3]{weih}.$\ \ _\square$

Clearly, the set $S^\dag$ corresponding to the set $S$ considered in the above proof belongs to a rather low subrecursive class. Possibly this can be used for obtaining some subrecursive refinement of Theorem \ref{semialg}, and a comparison of this refinement with the result from \mbox{\cite[Theorem 4.2]{TZ}} could be to the point then.

\section*{Acknowledgments}

Thanks are due to Ivan Georgiev who read a previous version of the paper and made several useful suggestions.

\vspace{7mm}
\noindent\footnotesize{Dimiter Skordev\\
Faculty of Mathematics and Informatics, Sofia University\\
blvd. James Bourchier 5, 1126 Sofia, Bulgaria\\
e-mail address:~ {\tt skordev@fmi.uni-sofia.bg}\\
home page:~ {\tt http://www.fmi.uni-sofia.bg/fmi/logic/skordev/}}  
\end{document}